\documentclass[oneside,reqno]{amsart}
\usepackage{amsmath}
\newtheorem{theorem}{Theorem}[section]
\newtheorem{lemma}[theorem]{Lemma}
\newtheorem{proposition}[theorem]{Proposition}

\theoremstyle{definition}
\newtheorem{definition}[theorem]{Definition}

\numberwithin{equation}{section}

\newcommand{\blankbox}[2]

\begin{document}
\title{Weighted $CBMO$ Estimates for Commutators of Matrix Hausdorff Operator on Heisenberg Group}
\thanks{
amnaarshad2@outlook.com (A. Ajaib); ahabbasi123@yahoo.com (A. Hussain)}
\author{Amna Ajaib $^1$, Amjad Hussain $^2$\\
$^{1,2}$Department of Mathematics, QUAID-E-AZAM University,
Islamabad, 310027, Pakistan}
\address{$^{1,2}$Department of Mathematics, QUAID-E-AZAM University,
Islamabad, 310027, Pakistan}




\begin{abstract}
In this article, we study the commutators of  Hausdorff operator and establish their boundedness on weighted Herz space in the setting of Heisenberg group.
\end{abstract}

\maketitle

\section{Introduction}
The matrix Hausdorff operators defined on $n$-dimensional Euclidean space $\mathbb{R}^n$ was firstly appeared in \cite{LL}, in the form:
\begin{equation}\label{1}
H_{\Phi,A}f(x)=\int_{\mathbb{R}^n}\Phi(y) f(xA(y))dy.
\end{equation}
 Taking into consideration the duality of $H^1$ and $BMO$, the authors in \cite{LL} shown that $H_{\Phi,A}$ is bounded on Hardy spaces. Subsequently, similar boundedness of $H_{\Phi,A}$ was reconsidered in \cite{L} using atomic decomposition of Hardy spaces. The above cited publications are important as their results are the first attempts to study the high dimensional Hausdorff operators on $H^1(\mathbb R^n).$ Recently, Liflyand and Miyachi \cite{LM0} extended these results on $H^{p}(\mathbb R^n)$ spaces with $0<p<1.$

 In $2012,$ Chen et al. \cite{CFL} modified the form of (\ref{1}) by replacing the kernel function $\Phi(y)$ with $\Phi(y)/|y|^n:$
\begin{equation}\label{2}
H_{\Phi,A}f(x)=\int_{\mathbb{R}^n}\frac{\Phi(y)}{|y|^n} f(xA(y))dy.
\end{equation}
As a sub case, when $A(y)=\rm{diag}[1/|y|,1/|y|,...,1/|y|],$ they give another definition of $n$-dimensional Hausdorff operator:
\begin{equation}\label{3}
H_{\Phi}f(x)=\int_{\mathbb{R}^n}\frac{\Phi(y)}{|y|^n} f(\frac{x}{|y|})dy.
\end{equation}
Their results includes the boundedness of Hausdorff operators on Hardy spaces, local Hardy spaces, Herz and Herz-type Hardy spaces with a conclusion that these operators has better performance on Herz-type Hardy spaces than their performance of Hardy spaces. In the same year, with different co-authors, Chen et al. \cite{CFZ} extended the problem of boundedness of $H_{\Phi, A}$ to the product of Hardy type spaces. The boundedness results regarding Hausdorff operators on $H^1(\mathbb R^n)$ were improved in \cite{CZ}. The continuity of (\ref{2}) on Morrey spaces, Hardy-Morrey spaces, Block spaces and rectangularly defined spaces has also been discussed in \cite{BL}, \cite{D}, \cite{H1} and \cite{CM}, respectively. Similarly, some results regarding the boundedness of $H_\Phi$ can be found in \cite{HG2,GQ,GW}.

In the same way, the study of commutators to integral operators is important as it has many applications in the theory of partial differential equations and in characterizing function spaces. An attempt has been made in \cite{HA2} to discuss the boundedness of commutators of $H_{\Phi,A},$ defined by:
\begin{equation}\label{4}
H_{\Phi,A}^bf(x)=\int_{\mathbb{R}^n}\frac{\Phi(y)}{|y|^n}(b(x)-b(xA(y))) f(xA(y))dy,
\end{equation}
on function spaces when the symbol function $b$ are either from Lipschitz space or central $BMO$ space. However, when the matrix $A(y)$ is diagonal, we get the commutators of $H_\Phi$ which were studied in \cite{HG1,W1} and \cite{HA}. For detailed history and other developments regarding Hausdorff operators, we refer the interested readers to the review articles \cite{CFW,L11}.

Besides the Euclidean space $\mathbb R^n,$ the matrix Hausdorff operator can be defined on $p$-adic linear space $\mathbb Q_p^n,$ which is a locally compact commutative group under addition (see, for instance, \cite{V1,V2}) and on the Heisenberg group $\mathbb H^n$ \cite{RFW,WF,M11}. Since, we are mainly concerned with the study of the commutators of Hausdorff operator defined on the Heisenberg group  $\mathbb H^n,$ therefore, it is mandatory to introduce this group briefly and the definition of matrix Hausdorff operator on it first.

With underlying manifold
$\mathbb{R}^{2n}\times\mathbb{R},$ the Heisenberg group $\mathbb{H}^{n}$
is the Lie group under the law of non-commutative multiplication
\begin{equation*}\begin{aligned}
  x\cdot  y&=(x_{1}, x_{2},..., x_{2n+1})\cdot(y_{1}, y_{2},...,y_{2n+1})\\
 &=\left(x_{1}+y_{1},...,x_{2n}+y_{2n},x_{2n+1}+y_{2n+1}+2\sum^{n}_{j=1}(y_{j}x_{n+j}-x_{j}y_{n+j})\right).
\end{aligned}\end{equation*}

The above definition suggests that for $ x\in \mathbb H^n,$ we have $ x\cdot 0= x$ and $ x\cdot- x= 0.$
Therefore, the identity and inverse elements of $\mathbb H^n$ are as in usual $\mathbb R^{2n+1}$ space. The basis for the corresponding Lie algebra is formed
by the vector fields
\begin{eqnarray*}
&X_{j}&=\frac{\partial} { \partial x_{j}} + 2x_{n+j} \frac{\partial}{\partial x_{2n+1}}, \quad{1\le j\le n},\\
&X_{n+j}&=\frac{\partial}{\partial x_{n+j}}-2x_{j}\frac{\partial}{\partial x_{2n+1}},\quad{1\le j\le n},\\
&X_{2n+1}&=\frac{\partial}{\partial x_{2n+1}}.\\
\end{eqnarray*}
The only non-vanishing commutator relations satisfied by these vector fields are
\begin{equation*}[X_{j},X_{n+j}]=-4X_{2n+1}, \quad{1\le j\le n}.
\end{equation*}

The dilation, on the Heisenberg group $\mathbb{H}^{n},$  is defined as
\begin{equation*}
\delta_{r}(x_{1},x_{2},...,x_{2n},x_{2n+1})=(rx_{1},rx_{2},rx_{2n},r^{2}x_{2n+1}),\quad r>0.
\end{equation*}
Also, the Haar measure on $\mathbb{H}^{n}$ coincides with the
usual Lebesgue measure on $\mathbb{R}^{2n}\times\mathbb{R}^1.$
Thus, for any measurable set $E\subset\mathbb{H}^{n},$ we denote its measure by $|E|.$
Moreover, it is easy to see that
\begin{equation*}
|\delta_{r}(E)|=r^{Q}|E|, \quad\quad d(\delta_r x)=r^{Q}dx,
\end{equation*}
where $Q=2n+2$ is the so-called homogeneous dimension of $\mathbb{H}^{n}.$

The Heisenberg group is a homogeneous group with the norm:
\begin{equation*}
| x|_{h}=\left[\left(\sum^{2n}_{i=1}x_i^{2}\right)^{2}+x^{2}_{2n+1}\right]^{1/4},
\end{equation*}
and the Heisenberg distance $d,$ generated by this norm is given by
\begin{equation*}
d( p, q)=d( q^{-1} p,0)=| q^{-1} p|_{h}.
\end{equation*}
Notice that $d$ satisfies triangular inequality and is left-invariant in the sense that
$$d( r\cdot  p, r\cdot  q)=d( p, q),\quad \forall ~ p, q, r\in\mathbb H^n.$$

The ball and sphere on $\mathbb H^n,$  for $r>0$ and $x\in\mathbb{H}^{n},$ can be defined as
\begin{equation*}
B(x,r)=\{y\in\mathbb{H}^{n}:d(x,y)<r\},
\end{equation*}
and
\begin{equation*}
S(x,r)=\{y\in\mathbb{H}^{n}:d(x,y)=r\},
\end{equation*}
respectively. To compute the measure of this ball on $\mathbb H^n,$ we proceed as below
\begin{equation*}
|B(x,r)|=|B(0,r)|=\Omega_{Q}r^{Q},
\end{equation*}
where $\Omega_{Q},$ being a function of $n$ only, is the volume of the unit ball $B(0,1).$ Also, the area of unit sphere $S(0,1)$ on
$\mathbb{H}^{n}$ is $w_{Q}=Q\Omega_{Q}.$
For further readings on Heisenberg group, we refer the interested reader to the book by Folland and Stein \cite {FS} and the works by authors in \cite{G,HPR,IMS}

Now, we are in position to define the Hausdorff operator and its commutators on the Heisenberg group $\mathbb{H}^{n}.$
Let $\Phi$ be locally integrable function on $\mathbb{H}^n$. The Hausdorff operators on $\mathbb{H}^n$
are defined by:
\begin{eqnarray*}\begin{aligned}[b]
T_{\Phi}f(x)&=&\int_{\mathbb{H}^n}\frac{\Phi(y)}{|y|^Q_h}f(\delta_{|y|^{-1}_h}x)dy,\\
T_{\Phi,A}f(x)&=&\int_{\mathbb{H}^n}\frac{\Phi(y)}{|y|^Q_h}f(A(y)x)dy,\\
\end{aligned}\end{eqnarray*}
where $A(y)$ is a matrix-valued function, and we assume that $\det A(y)\neq0$ almost every where in the support of $\Phi$.
Also, we define the commutators $T^{b}_{\Phi,A}$ of $T_{\Phi,A}$ with locally integrable function $b$ as
\begin{equation}
 T^{b}_{\Phi,A}(f)=b T_{\Phi,A}(f)-T_{\Phi,A}(bf).
\end{equation}

In this article, we will study the boundedness of $T^b_{\Phi,A}$ on the weighted Herz spaces $\dot{K}^{\alpha_2,p}_{q_2}(\mathbb H^n;w),$ defined in the next section, with the Heisenberg group as underlying space. The next section contains some basic definitions and notations likewise some necessary propositions which will be used in the succeeding sections. Finally, the last section is reserved for the main results of this study along with their proofs.
\section{\textbf{Some Definitions and Notations}}
In 1972, Muckenhoupt \cite{MH} studied the Hardy-Littlewood maximal function on weighted $L^p$ spaces and introduced the theory of $A_p$ weights as a result. The theory was well studied in the later work by Garc\'{\i}a-Cuerva et al. \cite{GRF}.
An extension of this theory, in the settings of Heisenberg group $\mathbb H^n,$ was provided in \cite{G} and studied in \cite{HPR,IMS}. Any non-negative, locally integrable function $w$ on $\mathbb H^n$ can be given the role of a weight. The notation $w(E)$ serves to define weighted measure of $E\subset\mathbb H^n,$ that is $w(E)=\int_Ew(x)dx.$
Also, if $p$ and $p^\prime$ satisfy $1/p+1/{p^\prime}=1,$ then they will be called mutually conjugate indices.

\begin{definition}\label{D1} We say that $w$ belongs to the Muckenhoupt class $A_p(\mathbb H^n),~1<p<\infty,$ if there exist a $C>0$ such that for
every ball $B\subset\mathbb H^n,$
\begin{eqnarray*}\left(\frac{1}{|B|}\int_Bw(x)dx\right)\left(\frac{1}{|B|}\int_Bw(x)^{-p^\prime/p}dx\right)^{p/p^\prime}\le C.\end{eqnarray*}
Also,  $w\in A_1$ if there exist a constant $C>0$ such that for
every ball $B\subset\mathbb H^n,$
$$\left(\frac{1}{|B|}\int_Bw(x)dx\right)\le C\operatorname*{essinf}_{x\in B} w(x).$$
When $p=\infty,$ we define $A_\infty=\bigcup_{1\le p<\infty}A_p.$
\end{definition}

According to Proposition 2.2 in \cite{RFW}, we have $A_p(\mathbb H^n)\subset A_q(\mathbb H^n),$ for $1\le p<q<\infty,$ and if $w\in A_p(\mathbb H^n),1<p<\infty,$ then there is an $\epsilon>0$ such that $p-\epsilon>1$ and $w\in A_{p-\epsilon}(\mathbb H^n).$ Therefore, we may use $q_w:=\inf\{q>1:w\in A_q\}$ to denote the critical index of $w.$

\begin{definition}\label{D11} We say that $w$ belongs to the
reverse H\"{o}lder class $RH_r(\mathbb H^n),$ if there exist a fixed constant $C>0$ and $r>1,$ such that for
every ball $B\subset\mathbb H^n,$
\begin{eqnarray*}\left(\frac{1}{|B|}\int_Bw^r(x)dx\right)^{1/r}\le\frac{C}{|B|}\int_Bw(x)dx.\end{eqnarray*}
\end{definition}
In \cite{IMS}, it was proved that $w\in A_\infty(\mathbb H^n)$ if and only if there exist some $r>1$ such that $w\in RH_r(\mathbb H^n).$
In addition, if $w\in RH_r(\mathbb H^n)$, $r>1$, then for some
$\epsilon>0$ we have $w\in RH_{r+\epsilon}(\mathbb H^n)$. We therefore use $r_w:= \sup\{r>1:w\in RH_r(\mathbb H^n)\}$ to denote the critical index of $w$ for the reverse H\"{o}lder condition.

A particular case of Muckenhoupt $A_p(\mathbb H^n)$ weights is the power weight function $|x|^\alpha_h$. From Proposition 2.3 in \cite{RFW}, for $x\in \mathbb H^n,$ we have $|x|^\alpha_h\in A_1(\mathbb H^n)$ if and only if $-Q<\alpha\le0.$ Also, for $1<p<\infty, |x|_h\in A_p(\mathbb H^n),$ if and only if $-Q<\alpha<Q(p-1).$
In view of these observation , it is easy to see that for $0<\alpha<\infty,$ $$|x|_h^\alpha\in \bigcap_{\frac{Q+\alpha}{Q}<p<\infty} A_p(\mathbb H^n),$$
 where $(Q+\alpha)/Q$ is known as the critical index of $|x|^\alpha_h.$

The following two Propositions, proved in \cite{RFW}, concerning $A_p(\mathbb H^n)$ weights will be useful in establishing weighted estimates for $T_{\Phi,A}^b$ on Herz-type spaces on $\mathbb H^n.$
\begin{proposition}\label{P1}
Let $w\in A_p\cap RH_r(\mathbb H^n),~p\ge 1$ and $r>1.$ Then there exist constants $C_1,C_2>0$ such that
$$C_1\left(\frac{|E|}{|B|}\right)^{p}\le\frac{w(E)}{w(B)}\le C_2\left(\frac{|E|}{|B|}\right)^{(r-1)/r},$$
for any measurable subset $E$ of a ball $B.$ In general, for any $\lambda>1,$
$$w(B(x_0,\lambda R)\le \lambda^{Qp}w(B(x_0,R)).$$
\end{proposition}
\begin{proposition}\label{P2}
 If $w\in A_p(\mathbb H^n),1\le p<\infty,$ then for any $f\in L^1_{\rm{loc}}(\mathbb H^n)$ and any ball $B\subset \mathbb H^n$
$$\frac{1}{|B|}\int_{B}|f(x)|dx\le C\left(\frac{1}{w(B)}\int_{B}|f(x)|^pw(x)dx\right)^{1/p}.$$
\end{proposition}

 For any measurable set $E\subset \mathbb H^n,$ the weighted Lebesgue space $L^p(E;w)$
 is the space of all functions $f$ satisfying the norm conditions
$$\|f\|_{L^p(E;w)}=\left(\int_E|f(x)|^pw(x)dx\right)^{1/p}<\infty,$$
where $1\le p<\infty$ and $w$ is a weight function  on $\mathbb H^n.$ When $p=\infty,$ we have
 $L^\infty(\mathbb H^n;w)=L^\infty(\mathbb H^n)$ and $\|f\|_{L^\infty(\mathbb H^n;w)}=\|f\|_{L^\infty(\mathbb H^n)}.$

Let $B_k:=\{x\in \mathbb H^n:|x|_h<2^k\}, E_k=B_k/B_{k-1}$ for $k\in \mathbb Z.$ Then the homogeneous weighted Herz space in the setting of Heisenberg group can be defined as below.
\begin{definition}
Let $\alpha\in\mathbb R, 0<p,q<\infty, $ and $w$ is a weight function on $\mathbb H^n.$ The homogeneous weighted Herz space $\dot{K}^{\alpha,p}_q(\mathbb H^n)$ is defined by
$$\dot{K}^{\alpha,p}_q(\mathbb H^n;w):=\left\{f \in L_{\rm loc}^q(\mathbb H^{n}/\{0\};w):\|f\|_{\dot{K}^{\alpha,p}_q(\mathbb H^n;w)}< \infty\right\},$$
where
$$\|f\|_{\dot{K}^{\alpha,p}_q(\mathbb H^n;w)}=\left\{\sum_{k=-\infty}^\infty w(B_k)^{\alpha p/Q}\|f\|_{L^p(E_k;w)}^p\right\}^{1/p}.$$
\end{definition}
When $w=1,$ we obtain $\dot{K}^{\alpha,p}_q(\mathbb H^n)$ introduced in \cite{LLU}. It is easy to verify that $\dot{K}^{\alpha,p}_p(\mathbb H^n)=L^p(\mathbb H^n,|\cdot|^{\alpha p}_h).$ Hence, Herz space can be considered as an extension of power weighted Lebesgue space.

\begin{definition}\label{D3}
Let $1<q<\infty$ and $w$ be a weight function on $\mathbb{H}^n$. Then, we say a function $f \in L^q_{loc}(\mathbb{H}^n;w)$ belongs to the weighted central bounded mean oscillation ($CBMO$) space $C\dot{M}O^q(\mathbb{H}^n;w)$ if
\begin{equation*}
\|f\|_{C\dot{M}O^q(\mathbb{H}^n;w)}=\underset{R>0}{\sup}\left(\frac{1}{w(B(0,R))}\int_{B(0,R)}|f(x)-f_B|^q w(x)\right)^{1/q}<\infty,
\end{equation*}
where
\begin{equation}
f_B=\frac{1}{|B(0,r)|}\int_{B(0,r)}f(x)dx.
\end{equation}
\end{definition}
\noindent For detailed study of $CBMO$ space on $\mathbb R^n,$ we refer the reader to \cite{AG,LY}.

Recently, weighted boundedness of matrix Hausdorff operators and their commutators defined on different underlying spaces are established in
\cite{RF,RF2,CHZ,HA1,SFH,RFW1,HS1,HS2}.
\begin{lemma}(\cite{RFW}) Suppose that the $(2n+1)\times(2n+1)$ matrix $M$ is invertible. Then
\begin{equation}\label{ES1}\|M\|^{-Q}\le|\det M^{-1}|\le\|M^{-1}\|^Q,\end{equation} where
\begin{equation}\|M\|=\underset{x\in\mathbb{H}^n,x\neq0}{\sup}\frac{|Mx|_h}{|x|_h}.\end{equation}
\end{lemma}

Also, when $A_p$ weights are reduced to the power function, we shall use the notation $v(\cdot)$ instead of $w(\cdot)$, that is $v(\cdot)=|\cdot|_h^\beta.$ In that case, an easy computation results in:
\begin{equation}\begin{aligned}
v(B_k)=\int_{|x|_h\leq2^k}\left|x\right|^{\beta}_h dx=\omega_Q 2^{k(Q+\beta)}/(\beta+Q).
\end{aligned}\end{equation}

Moreover, in case of boundedness of $T_{\Phi,A}^b$ on power weighted Herz space, we shall frequently use the piecewise defined function $G:$
\begin{eqnarray*}\begin{aligned}
 G(M,\delta\beta)&={\begin{cases} \|M\|^{\delta\beta}  & \text{if } \beta>0,\\
\|M^{-1}\|^{-\delta\beta}  & \text{if } \beta\le0,\end{cases}}
\end{aligned}\end{eqnarray*}
where $M$ is any invertible matrix, $\alpha\in\mathbb R$ and $\delta$ is a positive real number. Then it is easy to see that $G$ satisfies:
\begin{eqnarray}\begin{aligned}[b]\label{ES}G(M,\beta(1/q+1/p))=&G(M,\beta/q)G(M,\beta/p),
\end{aligned}\end{eqnarray}
where $p,q\in\mathbb{Z^+}.$

\begin{proposition}\label{P11} Suppose that the $(2n+1)\times(2n+1)$ matrix $M$ is invertible. Let $\beta>-n, \ v(x)=| x|_h^{\beta}$ and $ x\in \mathbb H^n,$ then
 \begin{eqnarray*}\begin{aligned}
 v(M  x)&\le{\begin{cases} \|M\|^\beta v({x}) & \text{if } \beta>0,\\
\|M^{-1}\|^{-\beta} v({x}) & \text{if } \beta\le0,\end{cases}}\\
&=G(M,\beta)v({x}).
\end{aligned}\end{eqnarray*}
\end{proposition}

From this point forward, the notations $A \preceq B$ will imply that $A \le C B,$ for some $C>0.$
 Similarly, for some positive constants $C_1$ and $C_2,$ if $A\le C_1B$ and $B\le C_2A,$ then we will write $A\simeq B.$ Also, we will use an obvious notation $\lambda B(0,R)=B(0,\lambda R),$ for $\lambda>0.$

\section{Main Results and Their Proofs} This section contains the main results of this study and the relevant proofs. Our first result is:
\begin{theorem}\label{T1}
Let $1\leq p,q,q_1,q_2\leq\infty$ and $\alpha_1, \alpha_2 \in \mathbb{R}$ with $\alpha_1<0.$ Suppose that $1/s=1/q_1+1/q$ and $\alpha_1/Q+1/q_1=\alpha_2/Q+1/q_2.$
In addition, let $w \in A_1$ with the critical index $r_w$ for the reverse H\"{o}lder condition and $s>q_2 r_w/(r_w-1).$\\
$(i)$ If $1/q_1+\alpha_1/Q\geq0$, then for any $1<\delta<r_w,$
\begin{equation*}\begin{aligned}[b]
\left\|T^b_{\Phi,A}f\right\|_{\dot{K}^{\alpha_2,p}_{q_2}(\mathbb{H}^n;w)}\leq K_1\left\|b\right\|_{C\dot{M}O^q(\mathbb{H}^n;w)}\left\|f\right\|_{\dot{K}^{\alpha_1,p}_{q_1}(\mathbb{H}^n;w)},
\end{aligned}\end{equation*}
where
\begin{eqnarray*}\begin{aligned}[b]
K_1&=\int_{\|A(y)\|<1}
\frac{|\Phi(y)|}{|y|^Q_h}\left(1+|\det A^{-1}(y)|^{1/q}\|A(y)\|^{Q/q}\right)\\
&\quad\times|\det A^{-1}(y)|^{1/q_1}\|A(y)\|^{-\alpha_1}\log\frac{2}{\|A(y)\|}dy\\
&+\int_{\|A(y)\|\ge1}\frac{|\Phi(y)|}{|y|^Q_h}\left(1+|\det A^{-1}(y)|^{1/q}\|A(y)\|^{Q/q}\right)\\
&\quad\times|\det A^{-1}(y)|^{1/q_1}\|A(y)\|^{Q/q_1-(\alpha_1+Q/q_1)(\delta-1)/\delta}\log2\|A(y)\|dy.
\end{aligned}\end{eqnarray*}
$(ii)$ If $\alpha_1/Q+1/q_1<0$, then for any $1<\delta<r_w$
\begin{equation*}\begin{aligned}[b]\label{E1}
\left\|T^b_{\Phi,A}\right\|_{\dot{K}^{\alpha_2,p}_{q_2}(\mathbb{H}^n;w)}\leq K_2\left\|b\right\|_{C\dot{M}O^q(\mathbb{H}^n;w)}\left\|f\right\|_{\dot{K}^{\alpha_1,p}_{q_1}(\mathbb{H}^n;w)}.
\end{aligned}\end{equation*}
where
\begin{eqnarray*}\begin{aligned}[b]
K_2&=\int_{\|A(y)\|\ge1}
\frac{|\Phi(y)|}{|y|^Q_h}\left(1+|\det A^{-1}(y)|^{1/q}\|A(y)\|^{Q/q}\right)\\
&\quad\times|\det A^{-1}(y)|^{1/q_1}\|A(y)\|^{-\alpha_1}\log2\|A(y)\|dy\\
&+\int_{\|A(y)\|<1}\frac{|\Phi(y)|}{|y|^Q_h}\left(1+|\det A^{-1}(y)|^{1/q}\|A(y)\|^{Q/q}\right)\\
&\quad\times|\det A^{-1}(y)|^{1/q_1}\|A(y)\|^{Q/q_1-(\alpha_1+Q/q_1)(\delta-1)/\delta}\log\frac{2}{\|A(y)\|}dy.
\end{aligned}\end{eqnarray*}
\end{theorem}
When general weights are reduced to power weights, then the next theorem is:
\begin{theorem}\label{T2}
Let $1\leq p<\infty, \ 1<q,q_1,q_2<\infty$ and $\beta>-n.$
If $1/q_2=1/q+1/q_1$ and $1/q+\alpha_2/Q=\alpha_1/Q,$ then we have

\begin{eqnarray*}\begin{aligned}[b]
&\|T^b_{\Phi,A}\|_{\dot{K}^{\alpha_2,p}_{q_2}(\mathbb{H}^n;v)}\preceq K_3\|b\|_{C\dot{M}O^q(\mathbb{H}^n;v)}
\|f\|_{\dot{K}^{\alpha_1,p}_{q_1}(\mathbb{H}^n;v)}.
\end{aligned}\end{eqnarray*}
where $K_3$ is
\begin{eqnarray*}\begin{aligned}[b]
K_3=\begin{cases}
\int_{\mathbb{H}^n}\Theta(y)\left(1+\log_2\left(\|A^{-1}(y)\|\|A(y)\|\right)\right)dy,&\text{if} \quad \alpha_1=0, \\
\int_{\mathbb{H}^n}\Theta(y)G\big(A^{-1}(y),\alpha_1(Q+\beta)/Q\big)dy,&\text{if}\quad \alpha_1\ne0.
\end{cases}
\end{aligned}\end{eqnarray*}
and
\begin{eqnarray*}\begin{aligned}[b]\Theta(y)&=\frac{|\Phi(y)|}{|y|^Q_h}|\det A^{-1}(y)|^{1/q_1}\Bigg(\log\frac{2}{\|A(y)\|}\chi_{\{\|A(y)\|<1\}}+\log2\|A(y)\|\chi_{\{\|A(y)\|\ge1\}}\Bigg)\\
&\quad\times G\big(A^{-1}(y),\beta/q_1\big)\bigg(1+ |\det A^{-1}(y)|^{1/q}G\big(A^{-1}(y),\beta/q\big)\|A(y)\|^{(Q+\beta)/q}\bigg).
\end{aligned}\end{eqnarray*}
\end{theorem}

\subsection{Proof of Theorem 1.2}
Here, we have to show that
\begin{eqnarray*}
\left\{\sum^{\infty}_{k=-\infty}w(B_k)^{\alpha_2p/Q}\left\|T^b_{\Phi,A}f\right\|^p_{L^{q_2}(E_k,w)}\right\}^{1/p}\preceq \left\|f\right\|_{\dot{K}^{\alpha_2,p}_{q_2}(\mathbb{H}^n;w)}.
\end{eqnarray*}
 By the Minkowski inequality and necessary splitting, the inner norm $\|T^b_{\Phi,A}f\|^p_{L^{q_2}(E_k,w)}$ can be approximated as:
\begin{eqnarray}\begin{aligned}[b]\label{E1}
&\left\|\left(H_{\Phi, A}^{b}f\right)\right\|_{L^{q_2}(E_k;w)}\\
&=\left\|\left(\int_{\mathbb{H}^{n}}\frac{\Phi(y)}{|y|_{h}^{Q}}(b(x)-b(A(y)x))f(A(y)x)dy\right)\right\|_{L^{q_2}(E_k;w)}\\
&\leq\int_{\mathbb{H}^{n}}\frac{\Phi(y)}{|y|_{h}^{Q}}\|(b(x)-b(A(y)x))f(A(y)x)\|_{L^{q_2}(E_k;w)}dy\\
&\leq\int_{\mathbb{H}^{n}}\frac{\Phi(y)}{|y|_{h}^{Q}}\|(b(x)-b_{B_k})f(A(y)x)\|_{L^{q_2}(E_k;w)}dy\\
&\quad+\int_{\mathbb{H}^{n}}\frac{\Phi(y)}{|y|_{h}^{Q}}\|(b(A(y)x)-b_{\|A(y)\|B_{k}})f(A(y)x)\|_{L^{q_2}(E_k;w)}dy\\
&\quad+\int_{\mathbb{H}^{n}}\frac{\Phi(y)}{|y|_{h}^{Q}}\|(b_{B_k}-b_{\|A(y)\|B_{k}})f(A(y)x)\|_{L^{q_2}(E_k;w)}dy\\
&=I_{1}+I_{2}+I_3.\end{aligned}
\end{eqnarray}
While targeting $I_1$, we first compute $\|(b(x)-b(A(y)x)) f(A(y)x)\|_{L^{q_2}(E_k;w)}.$
The condition $s>q_2 r_w/(r_w-1)$ implies that there exist $1<r<r_w$ such that $s=q_2 r^{\prime}.$ Therefore, by H\"{o}lder inequality and reverse H\"{o}lder condition, we have
\begin{eqnarray}\begin{aligned}[b]\label{E2}
&\|\left(b(\cdot)-b_{B_k}\right)f(A(y)\cdot)\|_{L^{q_2}(E_k;w)}\\
&=\left(\int_{E_k}\left|(b(x)-b_{B_k})f(A(y)x)\right|^{s}dx\right)^{1/s}\left(\int_{E_k}w(x)^rdx\right)^{1/rq_2}\\
&\preceq|B_k|^{-1/s}w(B_k)^{1/q_2}\|(b(\cdot)-b_{B_k})f(A(y)\cdot)\|_{L^{s}(E_k)}.\end{aligned}
\end{eqnarray}
Next, utilizing the condition $1/s=1/q+1/q_1$, we can have
\begin{eqnarray}
\begin{aligned}\label{E3}
\|(b(\cdot)-b_{B_k})f(A(y)\cdot)\|_{L^{s}(E_k)}\le \|b(\cdot)-b_{B_k}\|_{L^{q}(B_k)}\|f(A(y)\cdot)\|_{L^{q_1}(B_k)}.
\end{aligned}
\end{eqnarray}
In second factor, on the right side of the inequality (\ref{E3}), a change of variables along with Proposition \ref{P2} yields
\begin{eqnarray}\begin{aligned}[b]\label{E4}
\|f(A(y)\cdot)\|_{L^{q_1}(B_k)}
&=|\det A^{-1}(y)|^{1/q_1}\left(\int_{A(y)B_{k}}|f(x)|^{{q_1}}dx\right)^{1/q_1}\\
&\preceq|\det A^{-1}(y)|^{1/q_1}|B(0,2^k\|A(y)\|)|^{1/q_1}\\
&\quad\times\left(\frac{1}{w(B(0,2^k\|A(y)\|))}\int_{B(0,2^k\|A(y)\|)}|f(x)|^{q_1}w(x)dx\right)^{1/q_1}\\
&\preceq\left(|\det A^{-1}(y)|\|A(y)\|^{Q}|B_k|\right)^{1/q_1}\\
&\quad\times w(\|A(y)\|B_{k})^{-1/q_1}\|f\|_{L^{q_1}(\|A(y)\|B_{k};w)}.
\end{aligned}\end{eqnarray}
Similarly, the other factor on the right hand of the inequality (\ref{E3}), in view of Proposition \ref{P2}, gives
\begin{eqnarray}\label{E5}
\|b(\cdot)-b_{B_k}\|_{L^{q}(B_k)}\preceq&|B_k|^{1/q}\|b\|_{C\dot{M}O^q(\mathbb{H}^n; w)}.
\end{eqnarray}
Inequalities (\ref{E2}--\ref{E5}) together yield
\begin{eqnarray*}\begin{aligned}[b]
&\|(b(\cdot)-b_{B_k})f(A(y)\cdot)\|_{L^{q_2}(E_k; w)}\\
&\preceq \|b\|_{C\dot{M}O^q(\mathbb{H}^n; w)}\|f\|_{L^{q_1}(\|A(y)\|B_{k}; w)}\\
&\quad\times \left(|\det A^{-1}(y)|\|A(y)\|^{Q}\right)^{1/q_1}
\frac{w(B_k)^{1/q_2}}{w(\|A(y)\|B_{k})^{1/q_1}}.\end{aligned}
\end{eqnarray*}
Hence, estimation of $I_1,$ can be given as
\begin{eqnarray*}
\begin{aligned}[b]
I_1&\preceq\|b\|_{C\dot{M}O^q(\mathbb{H}^n; w)}\int_{\mathbb{H}^n}\frac{|\Phi(y)|}{|y|^Q_h}\left(|\det A^{-1}(y)|\|A(y)\|^{Q}\right)^{1/q_1}\\
&\qquad\times
\frac{w(B_k)^{1/q_2}}{w(\|A(y)\|B_{k})^{1/q_1}}\|f\|_{L^{q_1}(\|A(y)\|B_{k}; w)}dy
\end{aligned}
\end{eqnarray*}

Next, we fix to estimate $I_2,$ which is given by
\begin{eqnarray*}
I_2&=&\int_{\mathbb{H}^n}\frac{|\Phi(y)|}{|y|^Q_h}\left\|(b(A(y)\cdot)-b_{\|A(y)\|B_k})f(A(y)\cdot)\right\|_{L^{q_2}(E_k; w)}dy
\end{eqnarray*}
Since $s=q_2 r^{\prime},$ therefore we infer from (\ref{E2}) that
\begin{eqnarray}\begin{aligned}[b]\label{E6}
&\left\|\left(b(A(y)\cdot)-b_{\|A(y)\|B_k}\right)f(A(y)\cdot)\right\|_{L^{q_2}(E_k;w)}\\
&\preceq|B_k|^{-1/s}w(B_k)^{1/q_2}\|(b(A(y)\cdot)-b_{\|A(y)\|B_k})f(A(y)\cdot)\|_{L^{s}(E_k)}.\end{aligned}
\end{eqnarray}
Applying change of variables formula, Proposition (\ref{P2}) and H\"{o}lder's inequality, to have
\begin{eqnarray}\begin{aligned}[b]\label{E7}
&\left\|(b(A(y).)-b_{\|A(y)\|B_{k}})f(A(y).)\right\|_{L^{s}(E_k)}\\
&=|\det A^{-1}(y)|^{1/s}\left(\int_{A(y)B_{k}}|(b(x)-b_{\|A(y)\|B_{k}})f(x)|^{s}dx\right)^{1/s}\\
&\preceq|\det A^{-1}(y)|^{1/s}|\|A(y)\|B_{k}|^{1/s}\\
&\quad\times\left(\frac{1}{w(\|A(y)\|B_{k})}\int_{\|A(y)\|B_{k}}|(b(x)-b_{\|A(y)\|B_{k}})f(x)|^{s}w(x)dx\right)^{1/s}\\
&\preceq|\det A^{-1}(y)|^{1/s}|B_k|^{1/s}\|A(y)\|^{Q1/s}w(\|A(y)\|B_{k})^{-1/s}\\
&\quad\times\left(\int_{\|A(y)\|B_{k}}|b(x)-b_{\|A(y)\|B_{k}}|^{q}w(x)dx\right)^{1/q}\left(\int_{\|A(y)\|B_{k}}|f(x)|^{q_1}w(x)dx\right)^{1/q_1}\\
&\preceq|\det A^{-1}(y)|^{1/s}|B_k|^{1/s}\|A(y)\|^{Q1/s}w(\|A(y)\|B_{k})^{-1/q_1}\\
&\quad\times\|f\|_{L^{q_1}(\|A(y)\|B_{k};w)} \|b\|_{C\dot{M}O^q(\mathbb{H}^n;w)}.\end{aligned}
\end{eqnarray}
By virtue of (\ref{E6}) and (\ref{E7}), the expression for $I_2$ assumes the following form:
\begin{eqnarray*}
\begin{aligned}[b]
I_2&\preceq\left\|b\right\|_{C\dot{M}O^q(\mathbb{H}^n, w)}\int_{\mathbb{H}^n}\frac{|\Phi(y)|}{|y|^Q_h}\left(|\det A^{-1}(y)|\|A(y)\|^{Q}\right)^{1/s}\\
&\quad\times
\frac{w(B_k)^{1/q_2}}{w(\|A(y)\|B_{k})^{1/q_1}}\|f\|_{L^{q_1}(\|A(y)\|B_{k}; w)}dy.
\end{aligned}
\end{eqnarray*}

Now, the estimation of $I_3,$  given by
\begin{eqnarray*}
\begin{split}
I_3&=\int_{\mathbb{H}^n}\frac{|\Phi(y)|}{|y|^Q_h}\|f(A(y)\cdot)\|_{L^{q_2}(E_k)}|b_{B_k}-b_{\|A(y)\|B_{k}}|dy,
\end{split}
\end{eqnarray*}
requires the bounds for $\|f(A(y)\cdot)\|_{L^{q_2}(E_k)}$ and $|b_{B_k}-b_{\|A(y)\|B_{k}}|.$
First we consider $\|f(A(y)\cdot)\|_{L^{q_2}(E_k, w)}.$ In view of the condition $s=q_2 r^\prime,$ we use H\"{o}lder's inequality and reverse H\"{o}lder condition to obtain
\begin{eqnarray}\begin{aligned}[b]\label{E8}
\left\|f(A(y)\cdot)\right\|_{L^{q_2}(E_k, w)}
&\le\left(\int_{B_k}|f(A(y)x)|^{q_2}w(x)dx\right)^{1/q_2}\\
&\le\left(\int_{B_k}|f(A(y)x)|^{s}dx\right)^{1/s}\left(\int_{B_k}w(x)^{r}dx\right)^{1/rq_2}\\
&\preceq|B_k|^{-1/s}w(B_k)^{1/q_2}\|f(A(y)\cdot)\|_{L^{s}(B_k)}.
\end{aligned}\end{eqnarray}
Furthermore, the condition $1/s=1/q+1/q_1$ and the inequality (\ref{E4}) help us to write
\begin{eqnarray}\begin{aligned}[b]\label{E9}
\|f(A(y)\cdot)\|_{L^{s}(B_k)}&=|B_k|^{1/q}\|f(A(y)\cdot)\|_{L^{q_1}(B_k)}\\
&\preceq|B_k|^{1/s}\left(|\det A^{-1}(y)|\|A(y)\|^{Q}\right)^{1/q_1}\\
&\quad\times w(\|A(y)\|B_{k})^{-1/q_1}\|f\|_{L^{q_1}(\|A(y)\|B_{k};w)}.
\end{aligned}\end{eqnarray}
We combine the inequalities (\ref{E8}) and (\ref{E9}) to substitute the result in the expression for $I_3,$ which now becomes
\begin{eqnarray*}\begin{aligned}[b]
I_3&\preceq\int_{\mathbb{H}^n}\frac{|\Phi(y)|}{|y|^Q_h}\left(|\det A^{-1}(y)|\|A(y)\|^Q\right)^{1/q_1}\\
&\quad\times
\frac{w(B(0,2^k))^{1/q_2}}{w(\|A(y)\|B_{k}))^{1/q_1}}\left\|f\right\|_{L^{q_1}(w(\|A(y)\|B_{k}), w)}|b_{B_k}-b_{\|A(y)\|B_{k}}|dy.
\end{aligned}
\end{eqnarray*}

Now, it turns to bound $|b_{B_k}-b_{\|A(y)\|B_{k}}|.$ For this purpose, we split the integral as below:
\begin{eqnarray*}\begin{aligned}[b]
I_3&\preceq\int_{\|A(y)\|<1}|b_{B_k}-b_{\|A(y)\|B_{k}}|\Psi(y)dy+\int_{\|A(y)\|\geq1}|b_{B_k}-b_{\|A(y)\|B_{k}}|\Psi(y)dy\\
&=I_{31}+I_{32},\end{aligned}
\end{eqnarray*}
where, for the convenience's sake, we used the following notation:
$$\Psi(y)=\frac{|\Phi(y)|}{|y|^Q_h}\left(|\det A^{-1}(y)|\|A(y)\|^Q\right)^{1/q_1}
\frac{w(B(0,2^k))^{1/q_2}}{w(\|A(y)\|B_{k}))^{1/q_1}}\|f\|_{L^{q_1}(\|A(y)\|B_{k};w)}.$$

Further decomposition of integral for $I_{31}$ results in:
\begin{eqnarray*}\begin{aligned}
I_{31}=&\sum^{\infty}_{j=0}\int_{2^{-j-1}\leq\|A(y)\|<2^{-j}}\Psi(y)\\
&\quad\Big\{\sum^j_{i=1}|b_{2^{-i}B_k}-b_{2^{-i+1}B_k}|+|b_{2^{-j}B_k}-b_{\|A(y)\|B_{k}}|\Big\}dy.
\end{aligned}\end{eqnarray*}
The first term inside the curly brackets can be approximated using Proposition \ref{P2}, that is
\begin{eqnarray*}\begin{aligned}[b]
&|b_{2^{-i}B_k}-b_{2^{-i+1}B_k}|\\
&\leq\frac{1}{|2^{-i}B_k|}\int_{2^{-i}B_k}|b(y)-b_{2^{-i+1}B_k}|dy\\
&\leq\frac{1}{w(2^{-i}B_k)}\int_{2^{-i}B_k}|b(y)-b_{2^{-i+1}B_k}|w(y)dy\\
&\leq\frac{1}{w(2^{-i}B_k)}\left(\int_{2^{-i+1}B_k}|b(y)-b_{2^{-i+1}B_k}|^qw(y)dy\right)^{\frac{1}{q}}
\left(\int_{2^{-i+1}B_k}w(y)dy\right)^{1/q^{\prime}}\\
&\leq\frac{{w(2^{-i+1}B_k)}}{{w(2^{-i}B_k)}}\left(\frac{1}{w(2^{-i+1}B_k)}\int_{2^{-i+1}B_k}|b(y)-b_{2^{-i+1}B_k}|^q w(y)dy\right)^{\frac{1}{q}}\\
&\preceq\left\|b\right\|_{C\dot{M}O^q(\mathbb{H}^n;w)}.
\end{aligned}
\end{eqnarray*}

Similarly, for second term inside the curly brackets in the expression of $I_{31},$ we have
\begin{eqnarray*}\begin{aligned}[b]
|b_{2^{-j}B_k}-b_{\|A(y)\|B_{k}}|\preceq\|b\|_{C\dot{M}O^q(\mathbb{H}^n;w)}.
\end{aligned}\end{eqnarray*}

Therefore, we finish the estimation of $I_{31}$ by writing
\begin{eqnarray*}\begin{aligned}
I_{31}&\preceq\|b\|_{C\dot{M}O^q(\mathbb{H}^n;w)}\sum^{\infty}_{j=0}\int_{2^{-j-1}\leq\|A(y)\|<2^{-j}}\Psi(y)(j+1)dy\\
&\preceq\|b\|_{C\dot{M}O^q(\mathbb{H}^n;w)}\int_{\|A(y)\|<1}\Psi(y){\log\frac{2}{\|A(y)\|}}dy.
\end{aligned}
\end{eqnarray*}

In a similar fashion, the integral $I_{32}$ gives us
\begin{eqnarray*}\begin{aligned}[b]
I_{32}
&=\int_{\|A(y)\|\geq1}\Psi(y)|b_{B_k}-b_{\|A(y)\|B_{k}}|dy.\\
&=\sum^{\infty}_{j=0}\int_{2^{j}\leq\|A(y)\|<2^{j+1}}\Psi(y)
\Big\{\sum^j_{i=1}|b_{2^{i}B_k}-b_{2^{i+1}B_k}|+|b_{2^{j+1}B_k}-b_{\|A(y)\|B_{k}}|\Big\}dy\\
&\preceq\|b\|_{C\dot{M}O^q(\mathbb{H}^n;w)}\int_{\|A(y)\|\ge1}\Psi(y) \ {\log2\|A(y)\|}dy.
\end{aligned}\end{eqnarray*}

A combination of expressions for $I_1$, $I_2$, $I_{31}$ and $I_{32}$, gives
\begin{eqnarray*}\begin{aligned}[b]
&\left\|T^b_{\Phi,A}f\right\|_{L^{q_2}(E_k; w)}&\\&\preceq\left\|b\right\|_{C\dot{M}O^q(\mathbb{H}^n;w)}\\
&\quad\times\int_{\mathbb H^n}\frac{|\Phi(y)|}{|y|^Q_h}\left(|\det A^{-1}(y)|\|A(y)\|^Q\right)^{1/q_1}\left(1+|\det A^{-1}(y)|^{1/q}\|A(y)\|^{Q/q}\right)\\
&\quad\times\frac{w(B(0,2^k))^{1/q_2}}{w(\|A(y)\|B_{k}))^{1/q_1}}\|f\|_{L^{q_1}(\|A(y)\|B_{k};w)}\max\Big\{\log\frac{2}{\|A(y)\|}, \ \log(2\|A(y)\|)\Big\}dy.
\end{aligned}\end{eqnarray*}

Keeping in view the definition of Herz space, factors containing the index $k$ in the expression of $\Psi(y)$ are important. Therefore, to proceed further and to avoid repetition of unimportant factors relative to the Herz space, we have to modify and rename the expression for $\Psi.$ Hence, in the remaining of this paper we shall use the following notation:
\begin{eqnarray*}\begin{aligned}[b]\tilde{\Psi}(y)&=\frac{|\Phi(y)|}{|y|^Q_h}\left(|\det A^{-1}(y)|\|A(y)\|^Q\right)^{1/q_1}\\
&\quad\times\left(1+|\det A^{-1}(y)|^{1/q}\|A(y)\|^{Q/q}\right)\max\Big\{\log\frac{2}{\|A(y)\|}, \ \log2\|A(y)\|\Big\}.\end{aligned}\end{eqnarray*}
Then,
\begin{eqnarray*}\begin{aligned}[b]
&\left\|T^b_{\Phi,A}f\right\|_{L^{q_2}(E_k;w)}\\
&\leq\left\|b\right\|_{C\dot{M}O^{q}(\mathbb{H}^n;w)}
\int_{\mathbb H^n}\tilde{\Psi}(y)\frac{w(B_k)^{1/q_2}}{w(\|A(y)\|B_{k})^{1/q_1}}\|f\|_{L^{q_1}(\|A(y)\|B_{k};w)}dy.
\end{aligned}\end{eqnarray*}

Finally, we take into consideration the definition of Herz space and employ the Minkowski inequality to have
\begin{eqnarray}\begin{aligned}[b]\label{E18}
\left\|T^b_{\Phi,A}f\right\|_{\dot{K}^{\alpha_2,p}_{q_2}(\mathbb{H}^n;w)}&=\left\{\sum^{\infty}_{k=-\infty}w(B_k)^{\frac{\alpha_2 p}{Q}}\left\|T^b_{\Phi,A}f\right\|^p_{L^{q_2}(E_k; w)}\right\}^{1/p}\\
&\preceq\|b\|_{C\dot{M}O^q(\mathbb{H}^n;w)}\int_{\mathbb H^n}\tilde{\Psi}(y)\\
&\times\left\{\left(\sum^{\infty}_{k=-\infty}\frac{w(B_k)^{\alpha_2/Q+1/{q_2}}}{w(\|A(y)\|B_{k})^{1/q_1}}
\|f\|_{L^{q_1}(\|A(y)\|B_{k};w)}\right)^p\right\}^{1/p}dy.
\end{aligned}\end{eqnarray}
Comparing inequality (\ref{E18}) with the inequality (3.9) in \cite{RFW}, we found that the term inside the curly brackets is same in both of these inequalities, the only difference lies in the integrands outside the curly brackets along with a constant multiple $\|b\|_{C\dot{M}O^q(\mathbb{H}^n;w)}$ outside the integral. Therefore, the inequality (\ref{E18}) can be written as \cite{RFW}:
\begin{eqnarray}\begin{aligned}[b]\label{E19}
&\left\|T^b_{\Phi,A}f\right\|_{\dot{K}^{\alpha_2,p}_{q_2}(\mathbb{H}^n;w)}&\\
&\preceq\|b\|_{C\dot{M}O^q(\mathbb{H}^n;w)}\sum^{\infty}_{j=-\infty}\int_{2^{j-1}<\|A(y)\|\leq 2^j}\tilde{\Psi}(y)\\
&\quad\times\Bigg\{\sum^{\infty}_{k=-\infty}\Bigg[\left(\frac{w(B_k)}{w(B_{k+j})}\right)^{\alpha_1/Q+1/{q_1}}\\
&\quad\times
\sum^{j}_{l=-\infty}\Big(\frac{w(B_{k+j})}{w(B_{k+l})}\Big)^{\alpha_1/Q}w(B_{k+l})^{\alpha_1/Q}\|f\|_{L^{q_1}(E_{k+l};w)}\Bigg]^p\Bigg\}^{1/p}dy,
\end{aligned}\end{eqnarray}
where the condition $\alpha_1/Q+1/{q_1}=\alpha_2/Q+1/{q_2}$ is utilized in obtaining the last inequality.

Under the stated condition that $\alpha_1<0$ and $l\leq j,$ we use Proposition \ref{P1} to have
\begin{eqnarray}\begin{aligned}[b]\label{E20}
\Big(\frac{w(B_{k+j})}{w(B_{k+l})}\Big)^{\alpha_1/Q}\preceq\Big(\frac{|B_{k+j}|}{|B_{k+l}|}\Big)^{\alpha_1(\delta-1)/(Q\delta)}=
2^{(j-l)\alpha_1(\delta-1)/\delta}.
\end{aligned}\end{eqnarray}
for any $1< \delta < r_w $.

In view od Proposition \ref{P1}, if $\alpha_1/Q+1/q_1\geq 0,$ then
\begin{eqnarray}\begin{aligned}[b]\label{E21}\left(\frac{w(B_k)}{w(B_{k+j})}\right)^{\alpha_1/Q+1/q_1}\preceq\begin{cases}
2^{-jQ(\alpha_1/Q+1/q_1)},&\text{if} \quad j\le0, \\
2^{-jQ(\alpha_1/Q+1/q_1)(\delta-1)/\delta},&\text{if} \quad j>0,
\end{cases}\end{aligned}\end{eqnarray}
and if $\alpha_1/Q+1/q_1< 0,$ then
 \begin{eqnarray}\begin{aligned}[b]\label{E211}\left(\frac{w(B_k)}{w(B_{k+j})}\right)^{\alpha_1/Q+1/q_1}\preceq\begin{cases}
2^{jQ(\alpha_1/Q+1/q_1)(\delta-1)/\delta},&\text{if} \quad j\le0, \\
2^{jQ(\alpha_1/Q+1/q_1)},&\text{if} \quad j>0,
\end{cases}\end{aligned}\end{eqnarray}
for any $1< \delta < r_w .$

Thus, for $\alpha_1/Q+1/q_1\geq 0,$ from  inequalities (\ref{E19}--\ref{E21}),
for any $1< \delta < r_w $, we have
\begin{eqnarray*}\begin{aligned}[b]
&\left\|T^b_{\Phi,A}f\right\|_{\dot{K}^{\alpha_2,p}_{q_2}(\mathbb{H}^n;w)}&\\
&\preceq\|b\|_{C\dot{M}O^q(\mathbb{H}^n;w)}\sum^{0}_{j=-\infty}\int_{2^{j-1}<\|A(y)\|\leq 2^j}\tilde{\Psi}(y)\|A(y)\|^{-\alpha_1-Q/q_1}\\
&\quad\times\sum^j_{l=-\infty}2^{\alpha_1(j-l)(\delta-1)/\delta}
\Bigg\{\sum^{\infty}_{k=-\infty}
w(B_{k+l})^{\alpha_1p/Q}\|f\|^p_{L^{q_1}(E_{k+l};w)}\Bigg\}^{1/p}dy\\
&\quad+\|b\|_{C\dot{M}O^q(\mathbb{H}^n;w)}\sum^{\infty}_{j=1}\int_{2^{j-1}<\|A(y)\|\leq 2^j}\tilde{\Psi}(y)\|A(y)\|^{(\alpha_1+Q/q_1)(\delta-1)/\delta}\\
&\quad\times\sum^j_{l=-\infty}2^{\alpha_1(j-l)(\delta-1)/\delta}
\Bigg\{\sum^{\infty}_{k=-\infty}
w(B_{k+l})^{\alpha_1p/Q}\|f\|^p_{L^{q_1}(E_{k+l};w)}\Bigg\}^{1/p}dy.
\end{aligned}\end{eqnarray*}
Replacing $\tilde{\Psi}(y)$ with its value in the above inequality, we get
\begin{eqnarray*}\begin{aligned}[b]
&\left\|T^b_{\Phi,A}f\right\|_{\dot{K}^{\alpha_2,p}_{q_2}(\mathbb{H}^n;w)}&\\
&\preceq\|b\|_{C\dot{M}O^q(\mathbb{H}^n;w)}\|f\|_{\dot{K}^{\alpha_1,p}_{q_1}(\mathbb{H}^n;w)}\\
&\quad\times\Bigg\{\int_{\|A(y)\|<1}
\frac{|\Phi(y)|}{|y|^Q_h}\left(1+|\det A^{-1}(y)|^{1/q}\|A(y)\|^{Q/q}\right)\\
&\quad\times|\det A^{-1}(y)|^{1/q_1}\|A(y)\|^{-\alpha_1}\log\frac{2}{\|A(y)\|}dy\\
&\quad+\int_{\|A(y)\|\ge1}\frac{|\Phi(y)|}{|y|^Q_h}\left(1+|\det A^{-1}(y)|^{1/q}\|A(y)\|^{Q/q}\right)\\
&\quad\times|\det A^{-1}(y)|^{1/q_1}\|A(y)\|^{Q/q_1-(\alpha_1+Q/q_1)(\delta-1)/\delta}\log(2\|A(y)\|)dy\Bigg\}.
\end{aligned}\end{eqnarray*}
This completes the proof of Theorem \ref{T1}$(i).$

Similarly, when $\alpha_1/Q+1/q_1<0,$ by using inequalities (\ref{E19}), (\ref{E20}) and (\ref{E211}), the second part of Theorem \ref{T1}  can be proved easily. Hence, we complete the proof of Theorem \ref{T1}.

\subsection{Proof of Theorem 1.3}

Following the proof of Theorem \ref{T1}, we write:
\begin{eqnarray*}\begin{aligned}[b]
&\left\|T^b_{\Phi,A}\right\|_{L^{q_2}(E_k;v)}\leq J_1+J_2+J_3,
\end{aligned}\end{eqnarray*}
where $J_1, J_2, \ \text{and} \ J_3$ are as $I_1, I_2, \ \text{and} \ I_3$ in the previous theorem with $w(\cdot)$ is replaced by $v(\cdot)=|\cdot|^{\alpha}_h.$
Then by using the  H\"{o}lder inequality and change of variables, we obtain
\begin{eqnarray*}\begin{aligned}[b]
J_1&\le\int_{\mathbb{H}^n}\frac{|\Phi(y)|}{|y|^Q_h}\Big(\int_{E_k}|b(x)-b_{B_k}|^{q}v(x)dx\Big)^{1/q}
\Big(\int_{E_k}|f(A(y)x)|^{q_1}v(x)dx\Big)^{1/q_1}dy\\
&\leq \ v(B_k)^{1/q}\|b\|_{C\dot{M}O^q(\mathbb{H}^n;v)}\\
&\quad\times\int_{\mathbb{H}^n}\frac{|\Phi(y)|}{|y|^Q_h}|\det A^{-1}(y)|^{1/q_1}\Big(\int_{A(y)E_{k}}|f(z)|^{q_1}v(A^{-1}(y)z)dz\Big)^{1/q_1}dy.
\end{aligned}\end{eqnarray*}
Using Proposition \ref{P11}, we get
\begin{eqnarray*}\begin{aligned}[b]
J_1&\leq  v(B_k)^{1/q}\|b\|_{C\dot{M}O^q(\mathbb{H}^n;v)}\int_{\mathbb{H}^n}\frac{|\Phi(y)|}{|y|^Q_h}|\det A^{-1}(y)|^{1/q_1}\\
&\qquad\times\Bigg(\int_{A(y)E_{k}}|f(x)|^{q_1}G\big(A^{-1}(y),\beta/q_1\big)v(x)dx\Bigg)^{1/q_1}dy\\
&\le  v(B_k)^{1/q}\|b\|_{C\dot{M}O^q(\mathbb{H}^n;v)}\\
&\quad\times\int_{\mathbb{H}^n}\frac{|\Phi(y)|}{|y|^{Q}_h}|\det A^{-1}(y)|^{1/q_1}G\big(A^{-1}(y),\beta/q_1\big)\|f\|_{L^{q_1}(\|A(y)\|E_{k};v)}dy.
\end{aligned}\end{eqnarray*}

Next, the expression for $J_2$ can be written as:
\begin{eqnarray}\begin{aligned}[b]\label{E22}
J_2=\int_{\mathbb{H}^n}\frac{|\Phi(y)|}{|y|^Q_h}\left\|(b(A(y)\cdot)-b_{\|A(y)\|B_{k}})f(A(y)\cdot)\right\|_{L^{q_2}(E_k;v)}dy.
\end{aligned}\end{eqnarray}
Changing variables and using the condition $q_2/q+q_2/q_1=1,$ we get
\begin{eqnarray}\begin{aligned}[b]\label{LE}
&\Big\|(b(A(y)x)-b_{\|A(y)\|B_{k}})(f(A(y).))\Big\|_{L^{q_2}(E_k;v)}\\
&=\Bigg(\int_{E_k}\left|(b(A(y)x)-b_{\|A(y)\|B_{k}})f(A(y)x)\right|^{q_2}v(x)dx\Bigg)^{1/q_2}\\
&=|\det A^{-1}(y)|^{1/q_2}G\big(A^{-1}(y),\beta/q_2\big)\\
&\quad\times\Bigg(\int_{A(y)E_{k}}\left|(b(x)-b_{\|A(y)\|B_{k}})f(x)\right|^{q_2}v(x)dx\Bigg)^{1/q_2}\\
&\le|\det A^{-1}(y)|^{1/q_2}G\big(A^{-1}(y),\beta/q_2\big)\\
&\quad\times\Bigg(\int_{\|A(y)\|B_{k}}\left|b(x)-b_{\|A(y)\|B_{k}}\right|^{q}v(x)dx\Bigg)^{1/q}
\Bigg(\int_{A(y)E_{k}}\left|f(x)\right|^{q_1}v(x)dx\Bigg)^{1/q_1}\\
&=|\det A^{-1}(y)|^{1/q_2}G\big(A^{-1}(y),\beta/q_2\big)\\
&\quad\times v(\|A(y)\|B_k)^{1/q}\|b\|_{C\dot{M}O^q(\mathbb{H}^n;v)}\|f\|_{L^{q_1}(A(y)E_{k};v)}.
\end{aligned}\end{eqnarray}
It is easy to see that $v(\|A(y)\|B_k)=\|A(y)\|^{Q+\beta}v(B_k).$ Using property (\ref{ES}) and (\ref{LE}), the inequality (\ref{E22}) becomes:
\begin{eqnarray*}\begin{aligned}[b]
J_2&=v(B_k)^{1/q}\|b\|_{C\dot{M}O^q(\mathbb{H}^n;v)}\int_{\mathbb{H}^n}\frac{|\Phi(y)|}{|y|^Q_h}|\det A^{-1}(y)|^{1/q_2}\\
&\quad\times G\big(A^{-1}(y),\beta/q\big)G\big(A^{-1}(y),\beta/q_1\big) \|A(y)\|^{(Q+\beta)/q}\|f\|_{L^{q_1}(A(y)E_{k};v)}dy.
\end{aligned}\end{eqnarray*}

It remains to estimates $J_3$. A change of variables following the
 H\"{o}lder inequality and Proposition \ref{P11} gives us
\begin{eqnarray*}\begin{aligned}[b]
J_3&=\int_{\mathbb{H}^n}\frac{|\Phi(y)|}{|y|^Q_h}
\Big\|(b_{B_k}-b_{\|A(y)\|B_{k}})f(A(y)\cdot)\Big\|_{L^{q_2}(E_k;v)}dy\\
&=\int_{\mathbb{H}^n}\frac{|\Phi(y)|}{|y|^Q_h}\|f(A(y)x)\|_{L^{q_2}(E_k;v)}
\left|b_{B_k}-b_{\|A(y)\|B_{k}}\right|dy\\
&\quad\times\Bigg(\int_{A(y)E_{k}}|f(x)|^{q_2}v(x)dx\Bigg)^{1/q_2}\left|b_{B_k}-b_{\|A(y)\|B_{k}}\right|dy\\
&\le\int_{\mathbb{H}^n}\frac{|\Phi(y)|}{|y|^Q_h}|\det A^{-1}(y)|^{1/q_2}G\big(A^{-1}(y),\beta/q_2\big)\\
&\quad\times v(\|A(y)\|B_k)^{1/q}\|f\|_{L^{q_1}(A(y)E_{k};v)}\left|b_{B_k}-b_{\|A(y)\|B_{k}}\right|dy
\end{aligned}\end{eqnarray*}

Next, if $\|A(y)\|<1,$ then there exists an integer $j\geq0,$ such that
$$2^{-j-1}\le\|A(y)\|< 2^{-j}.$$
Therefore,
\begin{eqnarray*}\begin{aligned}\left|b_{B_k}-b_{\|A(y)\|B_{k}}\right|&\leq\sum_{i=1}^{j}|2^{-i}b_{B_{k}}-2^{-i+1}b_{B_{k}}|+|2^{-j}b_{B_{k}}-b_{A(y)B_{k}}|\\
&\le\|b\|_{C\dot{M}O^q(\mathbb{H}^n;v)}\log\frac{2}{\|A(y)\|}.
\end{aligned}\end{eqnarray*}
Similarly, for $\|A(y)\|\ge1,$ we have
\begin{eqnarray*}\begin{aligned}\left|b_{B_k}-b_{\|A(y)\|B_{k}}\right|
&\le\|b\|_{C\dot{M}O^q(\mathbb{H}^n;v)}\log2\|A(y)\|.
\end{aligned}\end{eqnarray*}

Hence
\begin{eqnarray*}\begin{aligned}[b]
J_3&\le v(B_k)^{1/q}\|b\|_{C\dot{M}O^q(\mathbb{H}^n;v)}\\
&\quad\times\int_{\mathbb{H}^n}\frac{|\Phi(y)|}{|y|^Q_h}|\det A^{-1}(y)|^{1/q_1}G\big(A^{-1}(y),\beta/q_2\big) G\big(A^{-1}(y),\beta/q\big)\|A(y)\|^{(Q+\beta)/q}\\
&\quad\times\Bigg(\log\frac{2}{\|A(y)\|}\chi_{\{\|A(y)\|<1\}}+\log2\|A(y)\|\chi_{\{\|A(y)\|\ge1\}}\Bigg) \|f\|_{L^{q_1}(A(y)E_{k};v)}dy.
\end{aligned}\end{eqnarray*}

Thus combining $J_1$, $J_2$ and $J_3,$ we get
\begin{eqnarray}\begin{aligned}[b]\label{E23}
&\|T^b_{\Phi,A}\|_{L^{q_2}(E_k;v)}\\&
\le v(B_k)^{1/q}\|b\|_{C\dot{M}O^q(\mathbb{H}^n;v)}\int_{\mathbb{H}^n}\frac{|\Phi(y)|}{|y|^Q_h}|\det A^{-1}(y)|^{1/q_1}\\
&\quad\times G\big(A^{-1}(y),\beta/q_1\big)\bigg(1+ |\det A^{-1}(y)|^{1/q}G\big(A^{-1}(y),\beta/q\big)\|A(y)\|^{(Q+\beta)/q}\bigg)\\
&\quad\times\Bigg(\log\frac{2}{\|A(y)\|}\chi_{\{\|A(y)\|<1\}}+\log2\|A(y)\|\chi_{\{\|A(y)\|\ge1\}}\Bigg) \|f\|_{L^{q_1}(A(y)E_{k};v)}dy.
\end{aligned}\end{eqnarray}

For the approximation of $\|f(\cdot) \|_{L^{q}(A(y)C_k)},$ we consider the method used in \cite{RFW}.
Hence, the definition of $E_k$ and (\ref{ES1}) implies that
\begin{eqnarray*}A(y)E_k\subset\{x:\|A^{-1}(y)\|^{-1}2^{k-1}<|x|_h<\|A(y)\|2^k\}.\end{eqnarray*}

Now, there exist an integer $l$ such that for any $y\in \rm{supp}(\Phi),$ we have
\begin{eqnarray}\label{E31}2^{l}<\|A^{-1}(y)\|^{-1}<2^{l+1}.\end{eqnarray}
Finally, the inequality $\|A^{-1}(y)\|^{-1}\le\|A(y)\|$ implies that there exist a non-negative integer $m$ satisfying:
\begin{eqnarray}\label{E32}2^{l+m}<\|A(y)\|<2^{l+m+1}.\end{eqnarray}

We infer from (\ref{E31}) and (\ref{E32}) that:
\begin{eqnarray*}\log_2(\|A(y)\|\|A^{-1}(y)\|/2)<m<\log_2(2\|A(y\|\|A^{-1}(y))\|).\end{eqnarray*}
Therefore,
\begin{eqnarray*}A(y)E_k\subset\{x:2^{l+k-1}<|x|<2^{k+l+m+1}\}.\end{eqnarray*}
Hence,
\begin{eqnarray}\label{EA}\begin{aligned}
\| f\|_{L^{q_2}(A(y)E_k;v)}\le\sum_{j=l}^{l+m+1}\| f\|_{L^{q_2}(E_{k+j};v)}.
\end{aligned}
\end{eqnarray}

Incorporating the inequality (\ref{EA}) into (\ref{E23}), we obtain
\begin{eqnarray}\begin{aligned}[b]\label{E24}
\|T^b_{\Phi,A}\|_{L^{q_2}(E_k;v)}
\le v(B_k)^{1/q}\|b\|_{C\dot{M}O^q(\mathbb{H}^n;v)}\int_{\mathbb{H}^n}{\Theta} (y)\sum_{j=l}^{l+m+1}\| f\|_{L^{q_2}(E_{k+j})}dy,
\end{aligned}\end{eqnarray}
where
\begin{eqnarray*}\begin{aligned}[b]\Theta(y)&=\frac{|\Phi(y)|}{|y|^Q_h}|\det A^{-1}(y)|^{1/q_1}\Bigg(\log\frac{2}{\|A(y)\|}\chi_{\{\|A(y)\|<1\}}+\log2\|A(y)\|\chi_{\{\|A(y)\|\ge1\}}\Bigg)\\
&\quad\times G\big(A^{-1}(y),\beta/q_1\big)\bigg(1+ |\det A^{-1}(y)|^{1/q}G\big(A^{-1}(y),\beta/q\big)\|A(y)\|^{(Q+\beta)/q}\bigg)
\end{aligned}\end{eqnarray*}
A use of Minkowski inequality and the condition $1/q+\alpha_2/Q=\alpha_1/Q,$ yield
\begin{eqnarray*}\begin{aligned}[b]
&\left\|T^b_{\Phi,A}\right\|_{\dot{K}^{\alpha_2,p}_{q_2}(\mathbb{H}^n,v)}\\
&\preceq\|b\|_{C\dot{M}O^q(\mathbb{H}^n;v)}\\
&\quad\times\Bigg\{\sum^{\infty}_{k=-\infty}\Bigg(v(B_k)^{1/q+\alpha_2/Q}\int_{\mathbb{H}^n}\Theta(y)\sum_{j=l}^{l+m+1}
\|f\|_{L^{q_1}(E_{k+j},v)}dy\Bigg)^p\Bigg\}^{1/p}\\
&\preceq\|b\|_{C\dot{M}O^q(\mathbb{H}^n;v)}\\
&\quad\times\int_{\mathbb{H}^n}\Theta(y)\sum_{j=l}^{l+m+1}v(B_{-j})^{\alpha_1/Q}\Bigg\{\sum^{\infty}_{k=-\infty}\Bigg(v(B_{k+j})^{\alpha_1/Q}
\|f\|_{L^{q_1}(E_{k+j},v)}\Bigg)^p\Bigg\}^{1/p}dy\\
&\preceq\|b\|_{C\dot{M}O^q(\mathbb{H}^n;v)}\|f\|_{\dot{K}^{\alpha_1,p}_{q_1}(\mathbb{H}^n;v)}\int_{\mathbb{H}^n}\Theta(y)\sum_{j=l}^{l+m+1}v(B_{-j})^{\alpha_1/Q}dy.
\end{aligned}\end{eqnarray*}

It is easy to see that
\begin{equation*}\begin{aligned}[b]
\sum_{j=l}^{l+m+1}v(B_{-j})^{\alpha_1/Q}\simeq\sum_{j=l}^{l+m+1}2^{-j\alpha_1(Q+\beta)/Q}.
\end{aligned}\end{equation*}
Next, for $\alpha_1=0$,
\begin{equation*}\begin{aligned}[b]
\sum_{j=l}^{l+m+1}2^{-j\alpha_1(Q+\beta)/Q}=m+2\preceq 1+\log_2\left(\|A^{-1}(y)\|\|A(y)\|\right),
\end{aligned}\end{equation*}
and for $\alpha_1\neq 0$,
\begin{eqnarray*}\begin{aligned}[b]
\sum_{j=l}^{l+m+1}2^{-j\alpha_1(Q+\beta)/Q}&\simeq2^{-l\alpha_1(Q+\beta)/Q}\\
&\preceq
\begin{cases}
\|A^{-1}(y)\|^{\alpha_1(Q+\beta)/Q},&\text{if} \quad \alpha_1>0, \\
\|A(y)\|^{-\alpha_1(Q+\beta)/Q},&\text{if}\quad \alpha_1<0,\\
\end{cases}\\
&=G\big(A^{-1}(y),\alpha_1(Q+\beta)/Q\big).
\end{aligned}\end{eqnarray*}

Therefore,
\begin{eqnarray*}\begin{aligned}[b]
&\|T^b_{\Phi,A}\|_{\dot{K}^{\alpha_2,p}_{q_2}(\mathbb{H}^n;v)}\\
&\preceq\|b\|_{C\dot{M}O^q(\mathbb{H}^n;v)}\|f\|_{\dot{K}^{\alpha_1,p}_{q_1}(\mathbb{H}^n;v)}\\
&\quad\times\begin{cases}
\int_{\mathbb{H}^n}\Theta(y)\left(1+\log_2\left(\|A^{-1}(y)\|\|A(y)\|\right)\right)dy,&\text{if} \quad \alpha_1=0, \\
\int_{\mathbb{H}^n}\Theta(y)G\big(A^{-1}(y),\alpha_1(Q+\beta)/Q\big)dy,&\text{if}\quad \alpha_1\ne0.
\end{cases}\\
&=K_3\|b\|_{C\dot{M}O^q(\mathbb{H}^n;v)}
\|f\|_{\dot{K}^{\alpha_1,p}_{q_1}(\mathbb{H}^n;v)}.
\end{aligned}\end{eqnarray*}
Thus we complete the proof of Theorem \ref{T2}.

\end{document}